\documentclass[12pt]{amsart}
 \usepackage[T1]{fontenc}
\usepackage{amscd}
\usepackage{amssymb}

\allowdisplaybreaks[1]

\newtheorem{thm}{Th\'eor\`eme}[section]
\newtheorem{rema}[thm]{Remarque}

\newtheorem{prop}[thm]{Proposition}

\newtheorem{defn}[thm]{D\'efinition}
\theoremstyle{definition}

\numberwithin{equation}{section}

\newcommand{\R}{\mathbb R}

\newcommand{\n}{\noindent}

\newcommand{\resumename}{R\'esum\'e}
\newenvironment{resume}{\narrower\footnotesize\bf
\noindent\resumename.\quad\footnotesize\rm}{\par\bigskip}

\font\hb=cmbx12
\newcommand\pt{\hbox{\hb .}}

\begin{document}

\title[Cohomologie des graphes ascendants]{Cohomologie de Chevalley des
graphes ascendants}
\author[W. Aloulou, D. Arnal et R. Chatbouri]{Walid Aloulou, Didier Arnal et
Ridha Chatbouri}

\begin{abstract}
The space $T_{poly}(\mathbb R^d)$ of all tensor fields on $\mathbb R^d$, equipped with the Schouten bracket is a Lie algebra. The subspace of ascending tensors is a Lie subalgebra of $T_{poly}(\mathbb R^d)$.

In this paper, we compute the cohomology of the adjoint representations of this algebra (in itself and $T_{poly}(\mathbb R^d)$), when we restrict ourselves to cochains defined by aerial Kontsevitch's graphs like in our previous work (Pacific J of Math, vol 229, no 2, (2007) 257-292). As in the  vectorial graphs case, the cohomology is freely generated by all the products of odd wheels.\\
\end{abstract}

\address{D\'epartement de Math\'ematiques, Institut Pr\'eparatoire aux Etudes
d'Ing\'enieurs de Sfax, Route Menzel Chaker Km 0.5, Sfax, 3018,
Tunisie} 
\email{Walid.Aloulou@ipeim.rnu.tn}
\address{
Institut de Math\'ematiques de Bourgogne\\
UMR CNRS 5584\\
Universit\'e de Bourgogne\\
U.F.R. Sciences et Techniques
B.P. 47870\\
F-21078 Dijon Cedex\\France} \email{Didier.Arnal@u-bourgogne.fr}
\address{
D\'epartement de Math\'ematiques\\
Unit\'e de Recherche Physique Math\'ematique\\
Fa\-cult\'e des Sciences de Monastir\\
Avenue de l'environnement\\
5019 Monastir\\
Tuni\-sie}
\email{Ridha.Chatbouri@ipeim.rnu.tn}

\keywords{Graphes de Kontsevich, Cohomologie de Chevalley}
\subjclass[2000]{17B56, 53D50, 05C90}

\thanks{Ce travail a \'et\'e effectu\'e dans le cadre des accords CMCU (programmes Hubert Curien Utique) 06/S 1502 et 09/G 1502.\\
D. Arnal remercie la Facult\'e des Sciences de Monastir, W. Aloulou et R. Chatbouri remercient l'Institut de Math\'ematiques de Bourgogne pour leurs chaleureuses hospitalit\'es.}


\maketitle

\begin{resume}
L'espace des tenseurs ascendants est une sous alg\`ebre de Lie de
l'alg\`ebre de Lie (gradu\'ee) $T_{poly}({\mathbb R}^d)$ des
champs de tenseurs sur ${\mathbb R}^d$ muni du crochet de
Schouten.

Dans cet article, on calcule la cohomologie des repr\'esentations
adjointes de cette sous alg\`ebre de Lie, en se restreignant \`a
des cocha\^\i nes d\'efinies par des graphes de Kontsevich
a\'eriens comme dans \cite {[AAC1]} et \cite {[AAC2]}. On retrouve
un r\'esultat analogue \`a celui de la cohomologie des graphes
vectoriels et linéaires: elle est librement engendrée par des
produits de roues de longueurs impaires.
\end{resume}


\

\section{Introduction}\label{sec1}

\

Notons $T_{poly}({\mathbb R}^d)$ l'espace des tenseurs
contravariants compl\`etement antisym\'etriques sur ${\mathbb
R}^d$. Cet espace, muni du crochet de Schouten et de la graduation
$deg(\alpha)=k-1$, si $\alpha$ est un $k$-tenseur, est une
alg\`ebre de Lie gradu\'ee.

Un problème naturel est le calcul de la cohomologie de Chevalley
de l'action adjointe de cette algèbre de Lie.

\

Pour aborder cette question, on se restreint à des cocha\^ines
$C_\Delta$ d\'efinies \`a partir des graphes de Kontsevich
a\'eriens $\Delta$. En effet, il semble probable qu'on définisse
ainsi la partie la plus fondamentale de la cohomologie cherchée.

M\^eme ainsi posée, la question reste très difficile. En
particulier, la cohomologie de l'action adjointe de la sous
algèbre des tenseurs homogènes, de la forme
$$\alpha=\sum_{1\leq i_1,\dots, i_p\leq d}\alpha^{i_1\dots i_p}(x)
\partial_{x_{i_1}}\wedge\dots\wedge\partial_{x_{i_p}}
$$
où les fonctions $\alpha^{i_1\dots i_p}(x)$ sont des polyn\^omes
homog\`enes de degr\'e $p$ sur ${\mathbb R}^d$, est non triviale.

\

D'un autre c\^oté, nous avons dans \cite {[AAC1]} et \cite
{[AAC2]} calculé la cohomologie des algèbres de Lie
$Vect(\mathbb{R}^d)$ des champs de vecteurs et
$\mathcal{L}_{poly}(\mathbb{R}^d)$ des tenseurs linéaires.

\

Cet article est un premier pas vers la localisation de la
cohomologie de $T_{poly}({\mathbb R}^d)$ à sa restriction aux
tenseurs homogènes.

Pour cel\`a, on définit une sous-alg\`ebre de Lie int\'eressante
de $T_{poly}({\mathbb R}^d)$, c'est celle des tenseurs ascendants:
un $p$-tenseur $\alpha$ est dit ascendant ou de type $(q\leq p)$,
avec $p,q\in\mathbb{N}$ si:
$$\alpha=\sum_{1\leq i_1,\dots, i_p\leq d}\alpha^{i_1\dots i_p}(x)
\partial_{x_{i_1}}\wedge\dots\wedge\partial_{x_{i_p}}
$$
où toute fonction $\alpha^{i_1\dots i_p}(x)$ est un polyn\^ome
homog\`ene de degr\'e $q$ sur ${\mathbb R}^d$ tel que
  $$\begin{cases}
    q<p & \text{si $p>1$}, \\
    q\leq p & \text{si $p\leq1$}.
  \end{cases}$$

\noindent Muni du crochet de Schouten, l'espace des tenseurs
ascendants sur $\R^d$ noté $T_{poly}^{(\pt\leq\pt)}(\R^d)$ est une
sous alg\`ebre de Lie gradu\'ee de $T_{poly}(\R^d)$ qui contient
$\mathfrak{gl}(\R^d)$.

 \

Comme dans \cite{[AAC1]} et \cite{[AAC2]}, on veut d\'eterminer la
cohomologie de Chevalley de l'action adjointe de cette sous
alg\`ebre dans $T_{poly}({\mathbb R}^d)$, en se restreignant \`a
des cocha\^ines $C_\Delta$ d\'efinies \`a partir des graphes de
Kontsevich a\'eriens $\Delta$. La restriction aux tenseurs
ascendants correspond \`a l'utilisation de graphes pour lesquels
chaque sommet $i$ est de type noté $(q\leq p)$, c'est à dire, tel
qu'exactement $p$ flèches sortent de $i$ et $q$ flèches
aboutissent sur $i$ avec $\begin{cases}
    q<p & \text{si $p>1$}, \\
    q\leq p & \text{si $p\leq1$}.
  \end{cases}$

 On appellera de tels graphes des graphes ascendants. La
cohomologie se traduit, alors, par un op\'erateur de cobord
$\partial$ d\'efini directement sur l'espace des graphes a\'eriens
ascendants. On calcule ici la cohomologie de cet op\'erateur
$\partial$.

\

La situation est donc similaire \`a celle des graphes vectoriels
et linéaires a\'eriens \'etudi\'es dans \cite{[AAC1]} et
\cite{[AAC2]}. Rappelons que pour construire la cocha\^ine
$C_\Delta$ dans le cas des graphes ascendants, on doit ajouter au
graphe a\'erien $\Delta$ autant de jambes que de pieds, suivant la
m\'ethode de \cite{[AGM]}, et ceci pour n'importe quel nombre de
pieds. On obtient ainsi pour chaque graphe a\'erien ascendant une
suite infinie de graphes non a\'eriens.

\

On calcule la cohomologie par les m\'ethodes de \cite{[AAC1]} et
\cite{[AAC2]}: apr\`es la d\'efinition de l'ordre et du symbole
d'une combinaison lin\'eaire sym\'etrique de graphes a\'eriens
ascendants, on d\'efinit une homotopie qui permet de ramener,
modulo des cobords, chaque cocycle \`a un graphe a\'erien formé
par des sommets simples, c'est à dire des sommets de type $(q\leq
p )$ avec $p\leq1$. Autrement dit, chaque sommet reçoit et fait
sortir au plus une flèche. La cohomologie de $\partial$
co\"{\i}ncide alors avec celle calcul\'ee dans \cite{[AAC1]} et
\cite{[AAC2]}: elle est librement engendr\'ee par des produits de
roues
$$
R_{2k_1+1}\times R_{2k_2+1}\times\dots\times R_{2k_p+1}
$$

\noindent de longueur $2k_i+1$ impaires et strictement croissantes
tels que $0\leq k_1<k_2<\dots<k_p$.

On en d\'eduit que la cohomologie de Chevalley de l'action
adjointe des tenseurs ascendants dans $T_{poly}({\mathbb R}^d)$
est donn\'ee par des cocha\^ines $C_\delta$ o\`u $\delta$ est une
combinaison de roues impaires de longueur strictement plus petite
que $2d$.

\


\section{Notations et d\'efinitions}\label{sec2}

\

On consid\`ere l'alg\`ebre de Lie gradu\'ee $T_{poly}(\R^d)$ des
champs de tenseurs contravariants compl\`etement antisym\'etriques
sur $\R^d$. Cette alg\`ebre est gradu\'ee par $deg$ (o\`u
$deg(\alpha)=k-1$, si $\alpha$ est un $k$-tenseur) et son crochet
est celui de Schouten qui s'\'ecrit:
$$
[\xi_1\wedge\dots\wedge\xi_k,\eta_1\wedge\dots\wedge\eta_\ell]_S=\sum_{i=1}^k
\sum_{j=1}^\ell(-1)^{k-i+j-1}\xi_1\wedge\cdots\widehat{\xi_i}\dots\wedge\xi_k
\wedge[\xi_i,\eta_j]\wedge\eta_1\wedge\cdots\widehat{\eta_j}\dots\wedge\eta_\ell.
$$
o\`u les $\xi_i$ et les $\eta_j$ sont des champs de vecteurs sur
$\R^d$.

Prenons deux tenseurs $\alpha=\displaystyle\sum_{i_1\dots
i_k}\alpha^{i_1\dots
i_k}\partial_{i_1}\wedge\dots\wedge\partial_{i_k}$ et
$\beta=\displaystyle\sum_{j_1\dots j_\ell}\beta^{j_1\dots
j_\ell}\partial_{j_1}\wedge\dots\wedge\partial_{j_\ell}$, on
d\'efinit l'op\'erateur $\nabla: T_{poly}(\R^d)\longrightarrow
Der(T_{poly}(\R^d))$ par

$$
\nabla_\alpha\beta=\sum_{r=1}^k~(-1)^{r-1}~\sum_{\begin{smallmatrix}i_1,\dots,i_k \\ j_1,\dots
,j_\ell\end{smallmatrix}}\alpha^{i_1\dots i_k}\left(\partial_{i_r}\beta^{j_1\dots
j_\ell}\right)
\partial_{i_1}\wedge\cdots\widehat{\partial_{i_r}}\dots\wedge\partial_{i_k}\wedge
\partial_{j_1}\wedge\dots\wedge\partial_{j_\ell}.
$$

Avec cet op\'erateur, on peut d\'efinir le crochet de Schouten en
posant simplement

$$
[\alpha,\beta]_S=(-1)^{deg(\alpha)}\nabla_\alpha\beta-(-1)^{(deg(\alpha)+1)
deg(\beta)}\nabla_\beta\alpha.
$$

$(T_{poly}(\R^d), [~~,~~]_S, deg)$ est une alg\`ebre de Lie
gradu\'ee: le crochet de Schouten est antisym\'etrique et il est
de degr\'e $0$.\\

Soient $(p,q)\in\mathbb{N}^2$, un $p$-tenseur $\alpha$ est dit de
type $(q\leq p)$ si et seulement si
$$\alpha=\sum_{1\leq i_1,\dots, i_p\leq d}\alpha^{i_1\dots i_p}(x)
\partial_{x_{i_1}}\wedge\dots\wedge\partial_{x_{i_p}}
$$
où chaque fonction $\alpha^{i_1\dots i_p}(x)$ est un polyn\^ome
homog\`ene sur ${\mathbb R}^d$ de degr\'e $q$ tel que
  $$\begin{cases}
    q<p & \text{si $p>1$}, \\
    q\leq p & \text{si $p\leq1$}.
  \end{cases}$$

Un tenseur de type $(q\leq p)$ est dit aussi un tenseur ascendant.

\noindent Muni du crochet de Schouten, l'espace
$T_{poly}^{(\pt\leq\pt)}(\R^d)$ des tenseurs ascendants est une
sous alg\`ebre de Lie gradu\'ee de $T_{poly}(\R^d)$.

\

Cependant, on peut voir chaque application $n$-lin\'eaire
antisym\'etrique comme une application $n$-lin\'eaire sym\'etrique
si on change la graduation. Posons $|\alpha|=deg(\alpha)+1$ et
$Q(\alpha,\beta)=(-1)^{deg(\alpha)}[\alpha,\beta]_S$, le crochet
de Schouten se transforme en une opération $Q$ qui est
$|~|$-symétrique. Dans ce cas, l'opérateur de cohomologie de
Chevalley s'écrit: (voir par exemple \cite{[AAC1]})

\begin{prop} {\rm (Cohomologie sym\'etris\'ee)}

\

Une application $n$-lin\'eaire $|~|$-sym\'etrique $C$ est une
$n$-cocha\^\i ne, son cobord $\partial C$ est donn\'e par
\begin{align*}
(\partial
C)(\alpha_0,\dots,\alpha_n)&=\sum_{i=0}^n\Big(\varepsilon_{|\alpha|}
(i,0\dots\hat{i}\dots n)(-1)^{|C|(|\alpha_i|-1)}\nabla_{\alpha_i}
C(\alpha_0,\dots\widehat{\alpha_i}\dots,\alpha_n)\cr
&+(-1)^{|C|}\varepsilon_{|\alpha|}(0\dots\hat{i}\dots n,i)
\nabla_{C(\alpha_0,\dots\widehat{\alpha_i}\dots,\alpha_n)}\alpha_i\Big)\cr
&-\sum_{i\neq
j}\varepsilon_{|\alpha|}(i,j,0,\dots\hat{i}\dots\hat{j} \dots
n)C\left(\nabla_{\alpha_i}\alpha_j,\alpha_0,\dots\widehat{\alpha_i}\dots
\widehat{\alpha_j}\dots,\alpha_n\right),
\end{align*}
où $\varepsilon_{|\alpha|} (i,0\dots\hat{i}\dots n)$ est la
signature de la permutation graduée
$(\begin{smallmatrix}\alpha_0,\dots,\alpha_n\\\alpha_i,\alpha_0,\dots\hat{i}\dots,
\alpha_n\end{smallmatrix})$ et plus précisement, le signe
$\varepsilon_{|\alpha|}$ sur la transposition
$(\alpha_i,\alpha_j)$ est:
$\varepsilon_{|\alpha|}(i,j)=(-1)^{|\alpha_i||\alpha_j|}\varepsilon_{|\alpha|}(j,i)$.

\end{prop}


\

\section{Graphes a\'eriens et cocha\^ines}\label{sec 3}

\

Les cocha\^\i nes \'etudi\'ees dans cet article sont construites
\`a partir de graphes comme dans \cite{[K]} et \cite{[AGM]}. Pour
des tenseurs de type $(\pt\leq\pt)$, on consid\`ere la classe
restreinte des graphes qui sont de type $(\pt\leq\pt)$.
Définissons la construction dans ce cas.\\

Un graphe $\Delta$ de type $(\pt\leq\pt)$ est un graphe a\'erien
de Kontsevich, d\'efini par la donn\'ee de $n$ sommets
num\'erot\'es $1,\dots,n$ et d'ar\^etes qui sont des fl\`eches
$\overrightarrow{ab}$ telles que chaque sommet $i$ est de type
$(f_i\leq\ell_i)$, c'est \`a dire, exactement $\ell_i$ fl\`eches
sortent de $i$ et $f_i$ fl\`eches aboutissent sur $i$ avec
  $\begin{cases}
  f_i<\ell_i & \text{si $\ell_i>1$}, \\
    f_i\leq \ell_i & \text{si $\ell_i\leq1$}.
  \end{cases}$.

 On appelle $Deb(i)$ la collection
des ar\^etes partant d'un sommet $i$, rang\'ees dans un ordre
donn\'e ($Deb(i)$ peut \^etre vide). Le graphe $\Delta$ est donc
une liste de $n$ listes de fl\`eches:
$$
\Delta=\left(\left(\overrightarrow{1a_1^1},\dots,\overrightarrow{1a^1_{\ell_1}}
\right),\dots,\left(\overrightarrow{na_1^n},\dots,\overrightarrow{na^n_{\ell_n}}
\right)\right).
$$

Pour construire la cocha\^\i ne, on ajoute des ``pieds'' (sommets
terrestres) et des ``jambes'' (flèches aboutissant à des pieds)
\`a $\Delta$. Soit $(m_1,\dots,m_n)\in\mathbb{N}^n$, on se donne
$|m|=m_1+\dots+m_n$ nouveaux sommets (pieds) not\'es
$b^1_1,\dots,b^1_{m_1},\dots,b^n_1,\dots,b^n_{m_n}$ et $|m|$
nouvelles fl\`eches (jambes) not\'ees
$\overrightarrow{1b^1_1},\dots,\overrightarrow{1b^1_{m_1}},\dots,\overrightarrow
{nb^n_1},\dots,\overrightarrow{nb^n_{m_n}}$. On d\'efinit le
graphe non a\'erien $\Gamma$ compatible avec la num\'erotation de
ses sommets, c'est \`a dire tel que les $\ell_1+m_1$ premi\`eres
fl\`eches partent du sommet 1, les $\ell_2+m_2$ suivantes du
sommet 2, etc...
$$
\Gamma=\left(\left(\overrightarrow{1a_1^1},\dots,\overrightarrow{1a^1_{\ell_1}},
\overrightarrow{1b^1_1},\dots,\overrightarrow{1b^1_{m_1}}\right),\dots,\left(
\overrightarrow{na_1^n},\dots,\overrightarrow{na^n_{\ell_n}},
\overrightarrow{nb^n_1},\dots,\overrightarrow{nb^n_{m_n}}\right)\right).
$$

Notons $(x^1,\dots,x^d)$ des points de $\R^d$. Si la fl\`eche
$\overrightarrow{e_\Gamma}$ de $\Gamma$ est \`a la $r^{i\grave
eme}$ place dans la liste $\Gamma$, on notera
$t_{\overrightarrow{e_\Gamma}}$ l'indice $t_r$. On notera donc
$x^{t_{\overrightarrow{e_\Gamma}}}$ pour $x^{t_r}$ et
$\partial_{t_{\overrightarrow{e_\Gamma}}}$ pour
$\frac{\partial}{\partial x^{t_r}}$.

Pour tout sommet $i$ de $\Delta$, la collection $Deb(i)$ sera la
suite $\left(\overrightarrow{ia_1^i}
,\dots,\overrightarrow{ia^i_{\ell_i}},\overrightarrow{ib^i_1},\dots,
\overrightarrow{ib^i_{m_i}}\right)$ des fl\`eches issues de $i$,
si $\alpha$ est un $\ell_i+m_i$-tenseur, on notera
$\alpha^{t_{Deb(i)}}$ sa composante:
$$
\alpha^{t_{Deb(i)}}=\alpha^{t_{\overrightarrow{ia^i_1}}\dots
t_{\overrightarrow{
ia^i_{\ell_i}}}t_{\overrightarrow{ib^i_1}}\dots
t_{\overrightarrow{ib^i_{m_i}}}}.
$$

De m\^eme, on note
$Fin(i)=(\overrightarrow{c_1^ii},\dots,\overrightarrow{c_{f_i}^ii})$
la suite des fl\`eches arrivant sur $i$ ($Fin(i)$ peut \^etre
vide). On définit pour chaque $i$, l'op\'erateur
$\partial_{Fin(i)}=\frac{\partial^{f_i}}{\partial_{t_{\overrightarrow
{c_1^ii}}}\dots\partial_{t_{\overrightarrow{c_{f_i}^ii}}}}$
 ( c'est l'op\'erateur
identit\'e si aucune fl\`eche n'arrive sur $i$).

Chaque graphe non a\'erien $\Gamma$ permet de d\'efinir une
application $B_\Gamma$ qui, pour toute famille
$\alpha_1,\dots,\alpha_n$ de tenseurs contravariants totalement
antisym\'etriques et d'ordres respectifs
$\ell_1+m_1,\dots,\ell_n+m_n$ sur $\mathbb{R}^d$, associe un
op\'erateur $|m|$-diff\'erentiel
$B_\Gamma(\alpha_1,\dots,\alpha_n)$ d\'efini par

$$
B_\Gamma(\alpha_1,\dots,\alpha_n)=\sum_{t_{\overrightarrow{e}}=1}^d
\prod_{i=1}^n\partial_{Fin(i)}\alpha_i^{t_{Deb(i)}}
\partial_{t_{\overrightarrow{1b^1_1}}}\wedge\dots\wedge
\partial_{t_{\overrightarrow{nb^n_{m_n}}}}.
$$

\n{\bf Exemple:}

\

 On prend l'exemple du graphe $\Gamma$ ci-dessous

$
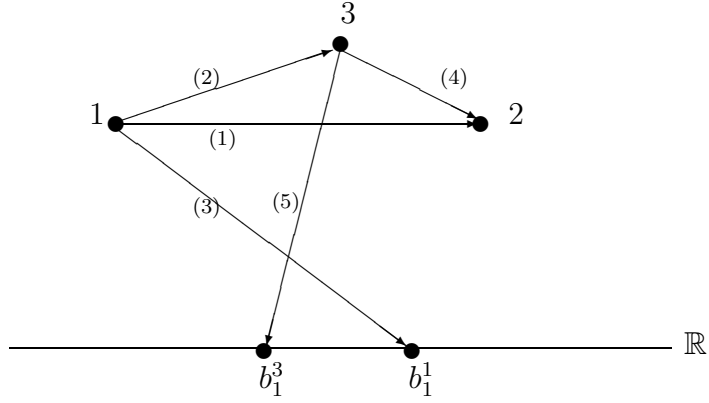
\begin{picture}(220,176)(-20,108)
\put(145,258){$3$}\put(60,220){\circle*{6}}\put(198,220){\circle*{6}}\put(145,250){\circle*{6}}
\put(89,237){$_{(2)}\quad\quad\quad\quad\quad\quad\quad_{(4)}$}
\put(58,220){\vector(3,1){84}}
\put(95,214){$_{(1)}$}\put(119,190){$_{(5)}$}
\put(89,188){$_{(3)}$}
\put(50,220){$1\quad\quad\quad\quad\quad\quad\quad\quad\quad\quad\quad\quad\quad
2$} \put(145,248){\vector(-1,-4){28}}
\put(145,248){\vector(2,-1){52}} \put(20,135){\line(1,0){250}}
\put(58,220){\vector(1,0){139}} \put(58,220){\vector(4,-3){112}}
\put(114,120){$b_1^3\quad\quad\quad\quad b_1^1$}
\put(116,134){\circle*{6}}\put(172,134){\circle*{6}}\put(275,133){$\mathbb{R}$}
\end{picture}$

Le graphe $\Gamma$ est orient\'e par l'ordre $\mathcal{O}$
suivant:
$$ \overrightarrow{12} < \overrightarrow{13}
< \overrightarrow{1b^1_1} < \overrightarrow{32} <
\overrightarrow{3b_1^3} .$$ Alors, pour tous tenseurs $\alpha_1,
\alpha_2$ et $\alpha_3$ d'ordre respectif $3,0$ et $2$,
l'op\'erateur bidiff\'erentiel $B_{\Gamma} $ s'\'ecrit:
$$B_{(\Gamma,\mathcal{O})} (\alpha_1 ,\alpha_2 ,
\alpha_3)  = \hskip-0.35cm\sum_{1\leq t_1,\dots,t_5 \leq
d}\hskip-0.25cm \alpha_1^{t_1,t_2,t_3} . \frac{\partial^2
\alpha_2}{\partial x^{t_1}\partial
x^{t_4}}.\frac{\partial\alpha_3^{t_4,t_5}}{\partial x^{t_2}}.
\frac{\partial }{\partial x^{t_3}} \wedge\frac{\partial}{\partial
x^{t_5} }.$$

\begin{rema}

\

La d\'efinition de l'op\'erateur $B_{\Gamma}$ d\'epend du choix de
la liste compatible

$\Gamma=\left(\overrightarrow{1a_1^1},\dots,\overrightarrow{1a^1_{\ell_1}},
\overrightarrow{1b^1_1},\dots,\overrightarrow{1b^1_{m_1}},\dots,
\overrightarrow{na_1^n},\dots,\overrightarrow{na^n_{\ell_n}},
\overrightarrow{nb^n_1},\dots,\overrightarrow{nb^n_{m_n}}\right)$.

\noindent Changer cette liste par une permutation $\sigma$ sur
l'ordre des fl\`eches faisant passer de la liste $\Gamma$ \`a la
liste $\sigma(\Gamma)$, revient \`a multiplier $B_{\Gamma}$ par le
signe $\varepsilon(\sigma)$.

\

La d\'efinition de l'op\'erateur $B_{\Gamma}$ est \'etendu aux
graphes non compatibles en posant
$$
B_{\sigma(\Gamma)}=\varepsilon(\sigma)B_{\Gamma}.
$$
\end{rema}

\

On considère un graphe aérien $
\Delta=\left(\left(\overrightarrow{1a_1^1},\dots,\overrightarrow{1a^1_{\ell_1}}
\right),\dots,\left(\overrightarrow{na_1^n},\dots,\overrightarrow{na^n_{\ell_n}}
\right)\right) $ de type $(\pt\leq\pt)$.
 Si on ajoute les jambes
\`a la fin de la liste $\Delta$, on obtient alors une liste
$$
\tilde{\Gamma}=\left(\overrightarrow{1a_1^1},\dots,\overrightarrow{1a^1_
{\ell_1}},\dots,\overrightarrow{na_1^n},\dots,\overrightarrow{n
a^n_{\ell_n}},\overrightarrow{1b^1_1},\dots,\overrightarrow{1b^1_{m_1}},
\dots,\overrightarrow{nb^n_1},\dots,\overrightarrow{nb^n_{m_n}}\right).
$$
On note $\varepsilon{(m_1,\dots,m_n)}$ la signature de la
permutation passant
de la liste $\Gamma$ \`a la liste $\tilde{\Gamma}$.\\

La cocha\^\i ne $C_\Delta$ d\'efinie par $\Delta$ est une s\'erie formelle:
$$
C_\Delta(\alpha_1,\dots,\alpha_n)=\sum_{m_1,\dots,m_n\geq0}\varepsilon{(m_1,
\dots,m_n)}B_{\tilde{\Gamma}}(\alpha_1,\dots,\alpha_n)
$$
avec la convention que $\alpha_i^{t_{Deb(i)}}=0$ si
$|\alpha_i|\neq\ell_i+m_i$.

Les cocha\^\i nes consid\'er\'ees ici sont des combinaisons
lin\'eaires \`a coefficients constants:
$$
C=\sum_{\Delta} a_\Delta C_\Delta
$$
de cocha\^\i nes associ\'ees \`a des graphes $\Delta$ de type
$(\pt\leq\pt)$. Comme les tenseurs $\alpha_i$ sont
antisym\'etriques, on se restreindra \`a des combinaisons $C$
telles que $a_{\sigma(\Delta)}=\varepsilon(\sigma)a_\Delta$ pour
toute permutation $\sigma$ sur l'ordre de l'une quelconque des
suites
$\left(\overrightarrow{ia_1^i},\dots,\overrightarrow{ia^i_{\ell_i}}\right)$.

Comme on ne consid\'ere que des cocha\^\i nes sym\'etriques (pour
$|~|$), on se restreint \`a des coha\^\i nes $C$ telles que pour
n'importe quelle permutation $\tau$ de deux sommets $i$ et $j$ de
$\Delta$, on a \
$a_{\tau(\Delta)}=(-1)^{\ell_i\ell_j+(\ell_i+\ell_j)(\ell_{i+1}+\dots+
\ell_{j-1})}a_\Delta$, si $\tau(\Delta)$ est le graphe ascendant:
\begin{align*}
&\tau(\Delta)=\cr&=\tau{\Big(}\Big(\overrightarrow{1a_1^1},\dots,
\overrightarrow{1a^1_{\ell_1}}\Big),\dots,\Big(\overrightarrow{ia^i_1},\dots,
\overrightarrow{ia^i_{\ell_i}}\Big),\dots,\Big(\overrightarrow{ja^j_1},\dots,
\overrightarrow{ja^j_{\ell_j}}\Big),\dots,\Big(\overrightarrow{na_1^n},\dots,
\overrightarrow{na^n_{\ell_n}}\Big){\Big)}\cr&={\Big(}\Big(
\overrightarrow{\tau(1)\tau(a_1^1)},\dots,\overrightarrow{\tau(1)\tau(a^1_{
\ell_1})}\Big),\dots,\Big(\overrightarrow{\tau(j)\tau(a^j_1)},\dots,
\overrightarrow{\tau(j)\tau(a^j_{\ell_j})}\Big),\dots\cr&\hskip3cm\dots,\Big(
\overrightarrow{\tau(i)\tau(a^i_1)},\dots,\overrightarrow{\tau(i)\tau
(a^i_{\ell_i})}\Big),\dots,\Big(\overrightarrow{\tau(n)\tau(a_1^n)},\dots,
\overrightarrow{\tau(n)\tau(a^n_{\ell_n})}\Big){\Big)}\cr&={\Big(}
\Big(\overrightarrow{1\tau(a_1^1)},\dots,\overrightarrow{1\tau(a^1_{\ell_1})}
\Big),\dots,\Big(\overrightarrow{i\tau(a^j_1)},\dots,\overrightarrow{i\tau(
a^j_{\ell_j})}\Big),\dots,\Big(\overrightarrow{j\tau(a^i_1)},\dots,
\overrightarrow{j\tau(a^i_{\ell_i})}\Big),\dots\cr&\hskip9.8cm\dots,\Big(
\overrightarrow{n\tau(a_1^n)},\dots,\overrightarrow{n\tau(a^n_{\ell_n})}\Big)
{\Big)}.
\end{align*}

Alors, d'apr\`es \cite{[AGM]}, on a
$$
C_{\tau(\Delta)}(\alpha_{1},\dots,\alpha_j,\dots,\alpha_i,\dots,\alpha_{n})=
C_{\Delta}(\alpha_1,\dots,\alpha_i,\dots,\alpha_j,\dots,\alpha_n).
$$

Si $\delta=\sum_\Delta a_\Delta\Delta$ est une combinaison
lin\'eaire sym\'etrique de graphes ascendants et $C_\delta$ la
cocha\^\i ne $\sum_\Delta a_\Delta C_\Delta$, alors, cette
condition garantit que le terme correspondant \`a
$m_1=\dots=m_n=0$ dans la s\'erie $C_\delta$ est
$|~|$-sym\'etrique. Gr\^ace au signe
$\varepsilon{(m_1,\dots,m_n)}$, tous les autres termes le sont
aussi.

\

\n{\bf Exemple:}

\

Prenons pour $\Delta$ la ``roue de longueur 3'':
$$
\Delta=\left(\left(\overrightarrow{12}\right),\left(\overrightarrow{23}\right),
\left(\overrightarrow{31}\right)\right).
$$
Alors si $\alpha_1,\alpha_2,\alpha_3$ sont des tenseurs d'ordre
$k_1,k_2,k_3$,
\begin{align*}
&C_\Delta(\alpha_1,\alpha_2,\alpha_3)=\cr&\sum_{1\leq t_1,\dots,
t_{k_1+k_2+k_3}\leq d
}\hskip-0.9cm(-1)^{k_2-1}\partial_{t_{k_1+k_2+1}}
\alpha_1^{t_1\dots t_{k_1}}\partial_{t_1}\alpha_2^{t_{k_1+1}\dots
t_{k_1+k_2}}
\partial_{t_{k_1+1}}\alpha_3^{t_{k_1+k_2+1}\dots t_{k_1+k_2+k_3}}\cr
&\hskip 1.9cm\partial_{t_2}
\wedge\dots\wedge\partial_{t_{k_1}}\wedge\partial_{t_{k_1+2}}\wedge\dots\wedge
\partial_{t_{k_1+k_2}}\wedge\partial_{t_{k_1+k_2+2}}\wedge\dots\wedge
\partial_{t_{k_1+k_2+k_3}}.
\end{align*}
Cette cocha\^\i ne n'est pas sym\'etrique, sa sym\'etrisation est:
$$
C(\alpha_1,\alpha_2,\alpha_3)=\sum_{\sigma\in
S_3}\varepsilon_{|\alpha|}(\sigma)
C_\Delta(\alpha_{\sigma(1)},\alpha_{\sigma(2)},\alpha_{\sigma(3)})=\left(3
C_\Delta-3C_{\Delta'}\right)(\alpha_1,\alpha_2,\alpha_3)$$ o\`u
$\Delta'$ est le graphe a\'erien:
$$
\Delta'=\left(\left(\overrightarrow{13}\right),
\left(\overrightarrow{21}\right),\left(\overrightarrow{32}\right)\right).
$$
$C$ est donc la cocha\^\i ne $C_\delta$ associ\'ee \`a la
combinaison sym\'etrique de graphes lin\'eaires
$\delta=3\Delta-3\Delta'$.


\section{L'op\'erateur de cobord sur les graphes de type $(\pt\leq\pt)$ }\label{sec4}

\

Dans ce paragraphe, nous d\'efinissons directement sur les
combinaisons lin\'eaires sym\'etriques $\delta$ de graphes de type
$(\pt\leq\pt)$, un op\'erateur de cobord not\'e $\partial$
correspondant \`a l'op\'erateur de cobord pour les cocha\^\i nes
sym\'etriques $C_\delta$ d\'efinies par $\delta$. Dans
\cite{[AGM]}, il est montr\'e qu'un tel op\'erateur
$\partial\delta$ d\'efini sur les graphes existe et v\'erifie
$C_{\partial\delta}=\partial C_\delta$. Dans le cas des graphes
ascendants, son expression peut \^etre simplifi\'ee.

\

Soit $\Delta$ un graphe de type $(\pt\leq\pt)$  ayant $n$ sommets
$(1,\dots,n)$.
$$
\Delta=\left(\left(\overrightarrow{1a_1^1},\dots,\overrightarrow{1a^1_{\ell_1}}
\right),\dots,\left(\overrightarrow{na_1^n},\dots,\overrightarrow{na^n_{\ell_n}}
\right)\right).
$$
Chaque sommet $k$ de $\Delta$ est de type $(f_k\leq \ell_k)$.

On fixe un indice $j$ tel que $0\leq j\leq n$. On renum\'erote les
sommets  de $\Delta$ en $(0,\dots\hat{j}\dots, n)$. On garde le
graphe avec son orientation, c'est \`a dire qu'on obtient un
graphe $\Delta(0,\dots\hat{j}\dots,n)$ qui a les m\^emes fl\`eches
que $\Delta$ mais leur nom a chang\'e conform\'ement au nouveau
nom des sommets. L'ordre des fl\`eches est inchang\'e.
$$
\Delta(0,\dots\hat{j}\dots,n)=\left(\left(\overrightarrow{0c_1^0},\dots,
\overrightarrow{0c^0_{\ell_1}}\right),\dots,\left(\overrightarrow{kc_1^k},\dots,
\overrightarrow{kc^k_{\ell_{k+1}}}\right),\dots,\hat{j},\dots,\left(
\overrightarrow{nc_1^n},\dots,\overrightarrow{nc^n_{\ell_n}}
\right)\right).
$$
Chaque sommet $k$ de $\Delta(0,\dots\hat{j}\dots,n)$ est de type
$(f_{k+1}\leq \ell_{k+1})$ si $k<j$ et $(f_k,\ell_k)$ si $k>j$.

\vskip 0.2cm

On fixe un autre indice $i<j$ et une partie $I$ de $Fin(i)$ dans
ce graphe $\Delta(0,\dots\hat{j}\dots,n)$. On construit les
graphes $\Delta_{i,j}^{I,r,s}$ et $\Delta_{j,i}^{I,r,t}$ non
compatibles de sommets $(0,\dots,n)$ ainsi:

\begin{align*}
&\Delta_{i,j}^{I,r,s}=\cr&\hskip-0.2cm\left(\left(\overrightarrow{0c_1^0}\dots
\overrightarrow{0c^0_{\ell_1}}\right),\dots,\left(\overrightarrow{ic_1^i}\dots
\underbrace{\overrightarrow{ij}}_{(s)}\dots\overrightarrow{ic^i_{r}}\right),
\left(\overrightarrow{jc_{r+1}^i}\dots\overrightarrow{jc^i_{\ell_{i+1}}}\right),
\dots\hat{j}\dots,\left(\overrightarrow{nc_1^n}\dots\overrightarrow{n
c^n_{\ell_n}}\right)\right)
\end{align*}
avec $I=Fin_{\Delta_{i,j}^{I,r,s}}(i)$,
$Fin_{\Delta_{i,j}^{I,r,s}}(j)=Fin_{\Delta(0,\dots\hat{j}\dots,n)}(i)\setminus
I$ et $s$ est la position de la fl\`eche $\overrightarrow{ij}$
dans $Deb_{\Delta_{i,j}^{I,r,s}}(i)$. On impose que $i$ soit de
type $(f'_i\leq r+1)$ et que $j$ soit de type
$(f'_j\leq\ell_{i+1}-r)$. On a $f'_i+f'_j=f_{i+1}+1$.

\

\noindent De m\^eme
\begin{align*}
&\Delta_{j,i}^{I,r,t}=\cr&\hskip-0.2cm\left(\left(\overrightarrow{0c_1^0}\dots
\overrightarrow{0c^0_{\ell_1}}\right),\dots,\left(\overrightarrow{jc_1^i}\dots
\underbrace{\overrightarrow{ji}}_{(t)}\dots\overrightarrow{jc^i_{r}}\right),
\left(\overrightarrow{ic_{r+1}^i}\dots\overrightarrow{ic^i_{\ell_{i+1}}}\right),
\dots\hat{j}\dots,\left(\overrightarrow{nc_1^n}\dots\overrightarrow{n
c^n_{\ell_n}}\right)\right)
\end{align*}
avec $I=Fin_{\Delta_{j,i}^{I,r,t}}(i)$,
$Fin_{\Delta_{j,i}^{I,r,t}}(j)=Fin_{\Delta(0,\dots\hat{j}\dots,n)}(i)\setminus
I$ et $t$ est la position de la fl\`eche $\overrightarrow{ji}$
dans $Deb_{\Delta_{j,i}^{I,r,t}}(i)$. On impose que $i$ soit de
type $(f'_i\leq r+1)$ et que $j$ soit de type
$(f'_j\leq\ell_{i+1}-r)$. On a $f'_i+f'_j=f_{i+1}+1$.
\\

Ces graphes sont d\'efinis pour tout $0\leq r\leq \ell_{i+1}$ et $1\leq
s,t\leq r+1$.

\

\

\begin{defn} {\rm (L'\'eclatement du sommet $i$)}

\

Soit $\Delta$ un graphe ayant $n$ sommets $(1,\dots,n)$. On fixe
deux indices $i$ et $j$ tels que $0\leq i<j\leq n$. On
renum\'erote les sommets  de $\Delta$ en $(0,\dots\hat{j}\dots,
n)$. On obtient un graphe $\Delta(0,\dots\hat{j}\dots, n)$
d\'efini comme pr\'ec\'edemment. L'\'eclatement du sommet $i$ de
$\Delta(0, \dots\hat{j}\dots,n)$ est la collection des graphes
$\Delta_{i,j}^{I,r,s}$ et $\Delta_{j,i}^{I,r,t}$, pour tout
$I\subset Fin_{\Delta(0,\dots\hat{j}\dots,n)}(i)$, $0\leq
r\leq\ell_{i+1}$ et $1\leq s,t\leq r+1$.

On dit que le graphe $\Delta_{i,j}^{I,r,s}$ (resp.
$\Delta_{j,i}^{I,r,t}$) se contracte proprement en $i$ \`a
$\Delta$ si et seulement si
$$
\inf\left\{\#Deb_{\Delta_{i,j}^{I,r,s}}(i)+\#Fin_{\Delta_{i,j}^{I,r,s}}(i),~\#Deb_{
\Delta_{i,j}^{I,r,s}}(j)+
\#Fin_{\Delta_{i,j}^{I,r,s}}(j)\right\}>1
$$
ou
$$
Deb_{\Delta_{i,j}^{I,r,s}}(i)=Fin_{\Delta_{i,j}^{I,r,s}}(j)=\left(\overrightarrow{ij}
\right)~\hbox{et}~Fin_{\Delta_{i,j}^{I,r,s}}(i)=Deb_{\Delta_{i,j}^{I,r,s}}(j)=\emptyset
$$
(resp.
$$
\inf\left\{\#Deb_{\Delta_{j,i}^{I,r,t}}(i)+\#Fin_{\Delta_{j,i}^{I,r,t}}(i),~
\#Deb_{\Delta_{j,i}^{I,r,t}}(j)+
\#Fin_{\Delta_{j,i}^{I,r,t}}(j)\right\}>1
$$
ou
$$
Deb_{\Delta_{j,i}^{I,r,t}}(j)=Fin_{\Delta_{j,i}^{I,r,t}}(i)=\left(\overrightarrow{ji}
\right)~\hbox{et}~Fin_{\Delta_{j,i}^{I,r,t}}(j)=Deb_{\Delta_{j,i}^{I,r,t}}(i)=
\emptyset).
$$

On notera la famille de graphes v\'erifiant cette propri\'et\'e:
$$
\{\Delta_{i,j}^{I,r,s}~\hbox{tels que
}~\Delta_{i,j}^{I,r,s}\begin{smallmatrix}prop\cr\rightarrow\cr
i\end{smallmatrix}
\Delta\}~~~~~~(\hbox{resp.}~\{\Delta_{j,i}^{I,r,t}~\hbox{tels que
}~\Delta_{j,i}^{I,r,t}\begin{smallmatrix}prop\cr\rightarrow \cr
i\end{smallmatrix} \Delta\}).
$$
\end{defn}

\

\begin{prop} {\rm (L'op\'erateur $\partial$ sur les graphes de type $(\pt\leq\pt)$)}

\

Soit $\delta=\sum_\Delta a_\Delta \Delta$ une combinaison
lin\'eaire sym\'etrique de graphes de type $(\pt\leq\pt)$ ayant
$n$ sommets $(1,\dots,n)$. Alors
$$
\partial C_\delta=\sum_\Delta a_\Delta C_{\partial\Delta}
$$
o\`u
\begin{align*}
\partial&(\Delta)=-\sum_{0\leq i<j\leq n}\left[\sum_{\Delta_{j,i}^{I,r,t}
\begin{smallmatrix}prop\cr\rightarrow\cr i\end{smallmatrix}\Delta}
\varepsilon_{j,i}^{r,t}.\Delta_{j,i}^{I,r,t}+\sum_{\Delta_{i,j}^{I,r,s}
\begin{smallmatrix}prop\cr\rightarrow\cr i\end{smallmatrix}\Delta}
\varepsilon_{i,j}^{r,s}.\Delta_{i,j}^{I,r,s}\right].
\end{align*}
Pour tout sommet $a$, posons $\#Deb_{\Delta_{j,i}^{I,r,t}}(a)=q_a$
et $\#Deb_{\Delta_{i,j}^{I,r,s}} (a)=p_a$. Les signes sont:
$$
\varepsilon_{j,i}^{r,t}=(-1)^{q_j(q_0+\dots+q_{j-1})}(-1)^{(q_j-1)(q_0+\dots+
q_{i-1})}(-1)^{t-1}
$$
et
$$
\varepsilon_{i,j}^{r,s}=(-1)^{p_j(p_0+\dots+p_{j-1})}(-1)^{(p_j-1)(p_0+\dots+
p_{i-1})}(-1)^{s-1}(-1)^{p_ip_j}.
$$

\end{prop}

\vskip 0.3cm
\noindent
{\bf Preuve}

\

La cohomologie de Chevalley est donn\'ee par:
\begin{align*}
\partial C_\delta(\alpha_0,\dots,\alpha_n)&=\sum_{j=0}^n
\varepsilon_{|\alpha|}(j,0\dots\hat{j}\dots
n)(-1)^{|C_\delta|(|\alpha_j|
-1)}\nabla_{\alpha_j}C_\delta(\alpha_0,\dots\widehat{\alpha_j}\dots,\alpha_n)\cr
&+\sum_{i=0}^n(-1)^{|C_\delta|}\varepsilon_{|\alpha|}(0\dots\hat{i}\dots
n,i)\nabla_{C_\delta(\alpha_0,\dots\widehat{\alpha_i}\dots,\alpha_n)}\alpha_i\cr
&-\sum_{i\neq
j}\varepsilon_{|\alpha|}(j,i,0,\dots\hat{i}\dots\hat{j} \dots
n)C_\delta\left(\nabla_{\alpha_j}\alpha_i,\alpha_0,\dots\widehat{\alpha_i}
\dots\widehat{\alpha_j}\dots,\alpha_n\right)\cr &=(I)+(II)+(III).
\end{align*}
Il y a des simplifications entre $(I)$, $(II)$ et $(III)$. Dans
\cite{[AGM]}, il est montr\'e que chaque terme de $(I)$ et $(II)$
correspond \`a un graphe $\Delta_{i,j}^{I,r,s}$ ou
$\Delta_{j,i}^{I,r,t}$ provenant d'un \'eclatement qui ne se
contracte pas proprement sur $\Delta$, on dira que ce terme n'est
pas propre. Il se simplifie avec un terme de $(III)$ et il nous
reste la somme des termes ``propres'' de $(III)$:
\begin{align*}
&\partial C_\delta(\alpha_0,\dots,\alpha_n)=-\sum_{\Delta}
a_{\Delta}
\left[\sum_{i<j}\varepsilon_{|\alpha|}(j,0,\dots\hat{j}\dots,n)
(-1)^{(|\alpha_j|-1)(|\alpha_0|+\dots+|\alpha_{i-1}|)}\right.\cr
&C_{\Delta}(\alpha_0,\dots,\underbrace{\nabla_{\alpha_j}
\alpha_i}_{(i)},
\dots\widehat{\alpha_j}\dots,\alpha_n)+\sum_{i<j}\varepsilon_{|\alpha|}(j,0,
\dots\hat{j}\dots,n)(-1)^{(|\alpha_j|-1)(|\alpha_0|+\dots+|\alpha_{i-1}|)}\cr
&(-1)^{|\alpha_i||\alpha_j|}\left.C_{\Delta}(\alpha_0,\dots,\underbrace{
\nabla_{\alpha_i}\alpha_j}_{(i)},\dots\widehat{\alpha_j}\dots,\alpha_n)\right].
\end{align*}
Apr\`es simplification, pour chaque graphe $\Delta$, l'op\'erateur
$$
C_{\Delta}(\alpha_0,\dots,\underbrace{\nabla_{\alpha_i}
\alpha_j}_{(i)}, \dots\widehat{\alpha_j}\dots,\alpha_n)
$$
est la somme des termes de la forme
$C_{\Delta_{i,j}^{I,r,s}}(\alpha_0,\dots, \alpha_n)$ o\`u
$\Delta_{i,j}^{I,r,s}$ se contracte proprement en $i$ sur
$\Delta$. De m\^eme l'op\'erateur
$$
C_{\Delta}(\alpha_0,\dots,\underbrace{\nabla_{\alpha_j}
\alpha_i}_{(i)}, \dots\widehat{\alpha_j}\dots,\alpha_n)
$$
est la somme des termes de la forme
$C_{\Delta_{j,i}^{I,r,t}}(\alpha_0,\dots, \alpha_n)$ o\`u
$\Delta_{j,i}^{I,r,t}$ se contracte proprement en $i$ sur
$\Delta$.

\

On v\'erifie que ces termes apparaissent avec les signes
$\varepsilon_{i,j}^{r,s}$ respectivement
$\varepsilon_{j,i}^{r,t}$.

\vskip 0.3cm
\begin{defn} {\rm (Espaces de cohomologie sur les graphes)}

\

Une combinaison lin\'eaire sym\'etrique $\delta$ de graphes de
type $(\pt\leq\pt)$ est un cocycle si $\partial\delta=0$, un
cobord s'il existe une combinaison lin\'eaire sym\'etrique $\beta$
telle que $\delta=\partial\beta$.

L'espace $Z^n$ des $n$-cocycles est l'espace des combinaisons
sym\'etriques de graphes lin\'eaires ayant $n$ sommets et qui sont
des cocycles.

L'espace $B^n$ des $n$-cobords est l'espace des combinaisons de
graphes qui sont les cobords de combinaisons sym\'etriques de
graphes ayant $n-1$ sommets.

Le $n^{\grave{e}me}$ espace de cohomologie des graphes $H^n$ est
le quotient de l'espace $Z^n$ par l'espace $B^n$.
\end{defn}


\

\section{Symbole d'un graphe}\label{sec5}

\

Soit $\Delta$ un graphe de type $(\pt\leq\pt)$ de sommets
num\'erot\'es
$(1,\dots,n)$. On peut distinguer cinq classes de sommets $i$:\\

Classe 1: les sommets $i$ de type $(0\leq\ell_i)$ avec $\ell_i>1$.
On
posera $\mathcal{O}(i)=(0,\ell_i)$.\\

Classe 2: les sommets $i$ de type $(f_i\leq\ell_i)$ avec
$\ell_i>f_i\geq1$. On
posera $\mathcal{O}(i)=(f_i,\ell_i)$.\\

Classe 3: les sommets $i$ de type $(1\leq1)$. On
posera $\mathcal{O}(i)=(1,1)$.\\

Classe 4: les sommets $i$ de type $(0\leq1)$. On
posera $\mathcal{O}(i)=(0,1)$.\\

Classe 5: les sommets $i$ de type $(0\leq0)$. On
posera $\mathcal{O}(i)=(0,0)$.\\

On ordonne les types des sommets suivant $(*)$ en posant:
\begin{align*}
&(0,\ell_i)>(f_j,\ell_j)>(1,1)>(0,1)>(0,0), \\& ~~(0,\ell_i)\geq
(0,\ell_{i'})\Leftrightarrow \ell_i\geq
\ell_{i'}~~\hbox{et}~~(f_i,\ell_i)\geq
(f_{i'},\ell_{i'})\Leftrightarrow
  \begin{cases}
    \ell_i\geq \ell_{i'} , \\
    \ell_i= \ell_{i'}~~ \hbox{et}~~ f_i\geq f_{i'}.
  \end{cases}
\end{align*}

L'ordre ${\mathcal O}(\Delta)$ d'un graphe  $\Delta$ de type
$(\pt\leq\pt)$ est le mot form\'e par les ordres de ses sommets:
$$
{\mathcal O}(\Delta)=\left({\mathcal O}(1),\dots,{\mathcal
O}(n)\right).
$$
On ordonne les ordres des graphes par l'ordre lexicographique, en
respectant $(*)$.

Si $\delta=\sum_{\Delta} a_{\Delta}\Delta$ est une combinaison
lin\'eaire sym\'etrique de graphes ascendants, on d\'efinit
l'ordre de $\delta$ par:
$$
{\mathcal O}(\delta)=Max\{{\mathcal O}(\Delta),~~a_{\Delta}\neq0\}
$$
et on appellera symbole de $\delta$ la combinaison lin\'eaire non
sym\'etrique:
$$
\sigma_\delta=\sum_{\begin{smallmatrix}\Delta\cr{\mathcal
O}(\Delta)={\mathcal O}(\delta)\end{smallmatrix}}a_{\Delta}\Delta.
$$
Comme $\delta$ est sym\'etrique, son ordre a la forme:
\begin{align*} &{\mathcal
O}(\delta)=\\&\big((0,\ell_1),...,(0,\ell_{k_0-1}),(f_{k_0},\ell_{k_0}),...,(f_{k_1-1},\ell_{k_1-1}),
\underset{(k_1)}{(1,1)},...,(1,1),\underset{(k_2)}{(0,1)},...,(0,1),
\underset{(k_3)}{(0,0)},...\big)
\end{align*} avec
$$
\ell_1\geq \ell_2\geq\dots\geq \ell_{k_0-1}, ~~\mbox{et}~~
~~(f_{k_0},\ell_{k_0})\geq\dots\geq (f_{k_1-1},\ell_{k_1-1})
$$
(on peut avoir $k_0=1$ ou $k_1=k_0$, etc\dots)\\

\vskip 0.3cm

\begin{prop} {\rm (Le symbole de $\partial\delta$)}

\

Soit $\delta$ une combinaison lin\'eaire sym\'etrique de graphes
a\'eriens ascendants, d'ordre
\begin{align*} &{\mathcal
O}(\delta)=\\&\big((0,\ell_1),...,(0,\ell_{k_0-1}),(f_{k_0},\ell_{k_0}),...,(f_{k_1-1},\ell_{k_1-1}),
\underset{(k_1)}{(1,1)},...,(1,1),\underset{(k_2)}{(0,1)},...,(0,1),
\underset{(k_3)}{(0,0)},...\big)
\end{align*}
Alors, chaque graphe $\Delta_{i,j}^{I,r,s}$ et
$\Delta_{j,i}^{I,r,t}$ apparaissant dans $\partial\delta$ est
d'ordre au plus:
\begin{align*}
&{\mathcal
O}(\delta)\oplus(1,1)=\cr&\big((0,\ell_1),...,(0,\ell_{k_0-1}),(f_{k_0},\ell_{k_0}),...,(f_{k_1-1},\ell_{k_1-1}),
\underset{(k_1)}{(1,1)},...,(1,1),\underset{(k_2+1)}{(0,1)},...,(0,1),
\underset{(k_3+1)}{(0,0)},...\big).
\end{align*}\\

Si ${\mathcal O}(\partial\delta)={\mathcal O}(\delta)\oplus(1,1)$,
alors le symbole de $\partial\delta$ est:
\begin{align*}
\sigma_{\partial\delta}&=-\sum_{\Delta\in\sigma_\delta}a_{\Delta}
\Big\{\sum_{\begin{smallmatrix}i<j\cr0\leq i<k_0-1\cr k_1-1\leq
j\leq k_2-1\end{smallmatrix}}\sum_{\Delta_{i,j}^{I,\ell_{i+1}-1,s}
\begin{smallmatrix}prop\cr\rightarrow\cr i\end{smallmatrix}\Delta}
\varepsilon_{i,j}^{\ell_{i+1}-1,s}\Delta_{i,j}^{I,\ell_{i+1}-1,s}\cr
&+\sum_{\begin{smallmatrix}i<j\cr k_0-1\leq i<k_1-1\cr k_1-1\leq
j\leq k_2-1\end{smallmatrix}}\Big[\sum_{\Delta_{j,i}^{I,0,1}
\begin{smallmatrix}prop\cr\rightarrow\cr i\end{smallmatrix}\Delta}
\varepsilon_{j,i}^{0,1}\Delta_{j,i}^{I,0,1}+\sum_{\Delta_{i,j}^{I,\ell_{i+1}-1,s}
\begin{smallmatrix}prop\cr\rightarrow\cr i\end{smallmatrix}\Delta}
\varepsilon_{i,j}^{\ell_{i+1}-1,s}\Delta_{i,j}^{I,\ell_{i+1}-1,s}\Big]\cr
&+\sum_{\begin{smallmatrix}i<j\cr k_1-1\leq i<k_2-1\cr k_1-1\leq
j\leq k_2-1\end{smallmatrix}}\Big[\sum_{\Delta_{j,i}^{0,1}
\begin{smallmatrix}prop\cr\rightarrow\cr i\end{smallmatrix}\Delta}
\varepsilon_{j,i}^{0,1}\Delta_{j,i}^{I,0,1}+\sum_{\Delta_{i,j}^{I,0,1}
\begin{smallmatrix}prop\cr\rightarrow\cr i\end{smallmatrix}\Delta}
\varepsilon_{i,j}^{0,1}
\Delta_{i,j}^{I,0,1}\Big]\Big\}\\&=\clubsuit.
\end{align*}
Réciproquement, si le second membre $\clubsuit$ est diff\'erent de
$0$, alors,
$\mathcal{O}(\partial\delta)=\mathcal{O}(\delta)\oplus(1,1)$ et il
sera \'egal \`a $\sigma_{\partial\delta}$.
\end{prop}

\vskip 0.3cm
\noindent
{\bf Preuve}

Fixons un couple $(i,j)$ avec $0\leq i<j\leq n$. Il est clair que
dans la d\'ecomposition de $\partial(\Delta)$, les ordres des
graphes $\Delta_{j,i}^{I,r,t}$ ou $\Delta_{i,j}^{I,r,s}$ sont:

{\bf *} Si $0\leq i<k_0-1$, alors

$$
{\mathcal
O}(\Delta_{j,i}^{I,r,t})=((0,\ell_1),\dots,\underbrace{(1,\ell_{i+1}-r)}_{(i)},\dots,
\underbrace{(0,r+1)}_{(j)},\dots)\qquad \hbox{où} ~~ 1\leq
r<\ell_{i+1}.
$$
Donc ${\mathcal O}(\Delta_{j,i}^{I,r,t})<{\mathcal
O}(\delta)\oplus(1,1)$. Il n'y a aucun graphe
$\Delta_{j,i}^{I,r,t}$, dans ce cas, qui appara\^\i t dans
$\clubsuit.$

D'autre part:
$$
{\mathcal
O}(\Delta_{i,j}^{I,r,s})=((0,\ell_1),\dots,\underbrace{(0,r+1)}_{(i)},\dots,
\underbrace{{(1,\ell_{i+1}-r)}}_{(j)},\dots)\qquad \hbox{où}
~~1\leq r<\ell_{i+1}.
$$
Donc ${\mathcal O}(\Delta_{i,j}^{I,r,s})\leq{\mathcal
O}(\delta)\oplus(1,1)$ et l'\'egalit\'e n'est vraie que si
$r=\ell_{i+1}-1$ et $k_1-1\leq j\leq k_2-1$. Pour chaque $s$
compris entre $1$ et $\ell_{i+1} $, il y a un seul graphe

\begin{align*}
&\Delta_{i,j}^{I,\ell_{i+1}-1,s}=\cr&\hskip-0.2cm\left(\left(\overrightarrow{0c_1^0}
\dots\overrightarrow{0c^0_{\ell_1}}\right),\dots,\left(\overrightarrow{ic_1^i}
\dots\underbrace{\overrightarrow{ij}}_{(s)}\dots\overrightarrow{i
c^i_{\ell_{i+1}-1}}\right),\left(\overrightarrow{jc^i_{\ell_{i+1}}}\right),\dots
\hat{j}\dots,\left(\overrightarrow{nc_1^n}\dots\overrightarrow{nc^n_{\ell_n}}
\right)\right)
\end{align*}

{\bf *} Si $k_0-1\leq i<k_1-1$, alors
$$
{\mathcal
O}(\Delta_{j,i}^{I,r,t})=((0,\ell_1),\dots,\underbrace{(f'_i,\ell_{i+1}-r)}_{(i)},\dots,
\underbrace{(f'_j,r+1)}_{(j)},\dots) \qquad \hbox{où} ~~0\leq r<
\ell_{i+1}.
$$
Donc ${\mathcal O}(\Delta_{j,i}^{I,r,t})\leq{\mathcal
O}(\delta)\oplus(1,1)$ et l'\'egalit\'e n'est vraie que si  $r=0$,
$t=1$, $f'_i=f_{i+1}$ et $k_1-1\leq j\leq k_2-1$. Il y a $f_{i+1}$
graphes $\Delta_{j,i}^{I,0,1}$ dans ce cas:
\begin{align*}
\Delta_{j,i}^{I,0,1}=\left(\left(\overrightarrow{0c_1^0}\dots\overrightarrow{0
c^0_{\ell_1}}\right),\dots,\left(\overrightarrow{ji}\right),\left(
\overrightarrow{ic_{1}^i}\dots\overrightarrow{ic^i_{\ell_{i+1}}}\right),\dots
\hat{j}\dots,\left(\overrightarrow{nc_1^n}\dots\overrightarrow{nc^n_{\ell_n}}
\right)\right)
\end{align*}
o\`u les $f_{i+1}$ fl\`eches arrivantes sur $i$ dans
$\Delta(0,\dots\hat{j}\dots,n)$ seront r\'eparties une sur $j$ et
les $f_{i+1}-1$ qui restent sur $i$ dans $\Delta_{j,i}^{I,0,1}$ et
on a $f_{i+1}$ possibilit\'es.

 D'autre part:
$$
{\mathcal
O}(\Delta_{i,j}^{I,r,s})=((0,\ell_1),\dots,\underbrace{(f'_i,r+1)}_{(i)},\dots,
\underbrace{{(f'_j,\ell_{i+1}-r)}}_{(j)},\dots)\qquad \hbox{où}
~~0\leq r<\ell_{i+1}.
$$
Donc ${\mathcal O}(\Delta_{i,j}^{I,r,s})\leq{\mathcal
O}(\delta)\oplus(1,1)$ et l'\'egalit\'e n'est vraie que si
$r=\ell_{i+1}-1$, $f'_i=f_{i+1}$ et $k_1-1\leq j\leq k_2-1$. Pour
chaque $s$ compris entre $1$ et $\ell_{i+1} $, il y a un seul
graphe
\begin{align*}
&\Delta_{i,j}^{I,\ell_{i+1}-1,s}=\cr&\hskip-0.2cm\left(\left(\overrightarrow{0c_1^0}
\dots\overrightarrow{0c^0_{\ell_1}}
\right),\dots,\left(\overrightarrow{ic_1^i}\dots
\underbrace{\overrightarrow{ij}}_{(s)}\dots\overrightarrow{ic^i_{\ell_{i+1}-1}}
\right),\left(\overrightarrow{jc^i_{\ell_{i+1}}}
\right),\dots\hat{j}\dots,\left(\overrightarrow{nc_1^n}\dots\overrightarrow{n
c^n_{\ell_n}}\right)\right)
\end{align*}

{\bf *} Si $k_1-1\leq i<k_2-1$, alors
$$
{\mathcal O}(\Delta_{j,i}^{I,r,t})={\mathcal
O}(\Delta_{i,j}^{I,r,s})=((0,\ell_1),\dots,
\underbrace{(1,1)}_{(i)},\dots,\underbrace{(1,1)}_{(j)},\dots)\qquad
\hbox{où} ~~r=0 ~~\mbox{et} ~~ s=t=1.
$$
Donc pour et seulement pour $i<j\leq k_2-1$, ${\mathcal
O}(\Delta_{j,i}^{I,0,1})= {\mathcal
O}(\Delta_{i,j}^{I,0,1})={\mathcal O}(\delta)\oplus(1,1)$ et il y
a un seul graphe $\Delta_{ji}^{0,1}$ et un seul graphe
$\Delta_{i,j}^{I,0,1}$ dans ce cas.
\begin{align*}
\Delta_{j,i}^{I,0,1}=\left(\left(\overrightarrow{0c_1^0}\dots\overrightarrow{0
c^0_{\ell_1}}\right),\dots,\left(\overrightarrow{ji}\right),\left(
\overrightarrow{ic^i_{1}}\right),\dots\hat{j}\dots,\left(\overrightarrow{nc_1^n}
\dots\overrightarrow{nc^n_{\ell_n}}\right)\right)
\end{align*}

\begin{align*}
\Delta_{i,j}^{I,0,1}=\left(\left(\overrightarrow{0c_1^0}\dots\overrightarrow{
0c^0_{\ell_1}}\right),\dots,\left(\overrightarrow{ij}\right),\left(
\overrightarrow{jc^i_{1}}\right),\dots\hat{j}\dots,\left(\overrightarrow{nc_1^n}
\dots\overrightarrow{nc^n_{\ell_n}}\right)\right)
\end{align*}

{\bf *} Si $i\geq k_2$, pour tout $j>i$, on a ${\mathcal
O}(\Delta_{j,i}^{I,r,t}), {\mathcal O}(\Delta_{i,j}^{I,r,s})
<{\mathcal O}(\delta)\oplus(1,1)$, dans ce cas, ces graphes
n'apparaissent pas dans $\clubsuit$.

\

Enfin, si $\Delta$ est un graphe de $\delta$ qui n'appara\^\i t
pas dans le symbole $\sigma_\delta$, alors les ordres des graphes
$\Delta_{i,j}^{I,r,s}$ et $\Delta_{j,i}^{I,r,t}$ qui se
contractent sur $\Delta$ sont tous strictement plus petits que
${\mathcal O}(\delta)\oplus(1,1)$. On regarde donc seulement les
$\Delta$ qui apparaissent dans $\clubsuit$.\vskip 0.2cm

A partir de maintenant, on note abusivement
$$
\sigma_{\partial\delta}=\sum_{\begin{smallmatrix}\Delta'\in\partial\delta\cr
{\mathcal O}(\Delta')={\mathcal
O}(\delta)\oplus(1,1)\end{smallmatrix}}a_{\Delta'} \Delta' ~~,~~
\hbox{si}~~\partial\delta=\sum_{\Delta'}a_{\Delta'}\Delta'.$$

\vskip 0.4cm

D\'efinissons maintenant l'op\'erateur d'homotopie.
\vskip 0.24cm

\begin{defn} {\rm (L'homotopie)}
\vskip 0.3cm

Soit $\Delta$ un graphe de type $(\pt\leq\pt)$ de sommets
$(0,\dots,n)$, ayant des sommets $i$ d'ordre $(1,1)$. On d\'efinit
l'homotopie $h$ comme l'application qui transforme le graphe
$\Delta$ en le graphe $h(\Delta)$ d\'efini ainsi:

On consid\`ere $i_0$, le plus grand des indices $i$ d'ordre
$(1,1)$. La fl\`eche issue de $i_0$ est $\overrightarrow{i_0a}$.
Le graphe $h(\Delta)$ est le graphe de sommets
$(0,\dots\hat{i_0}\dots,n)$ obtenu en comprimant la fl\`eche
$\overrightarrow{i_0a}$ de $\Delta$ et en identifiant les sommets
$i_0$ et $a$ au sommet $a$. On remplace la fl\`eche
$\overrightarrow{ci_0}$ dans $\Delta$  par la fl\`eche
$\overrightarrow{ca}$ dans $h(\Delta)$.

Si $\Delta$ n'a pas de sommets d'ordre $(1,1)$, on pose
$h(\Delta)=0$.

On prolonge $h$ lin\'eairement \`a l'espace des combinaisons lin\'eaires de
graphes.
\end{defn}

\vskip 0.3cm
\begin{prop} {\rm (Symbole et homotopie)}

\

Soit $\delta$ une combinaison lin\'eaire sym\'etrique de graphes
de type $(\pt\leq\pt)$ ayant $n$ sommets num\'erot\'es
$(1,\dots,n)$ et d'ordre
\begin{align*} &{\mathcal
O}(\delta)=\\&\big((0,\ell_1),...,(0,\ell_{k_0-1}),(f_{k_0},\ell_{k_0}),...,(f_{k_1-1},\ell_{k_1-1}),
\underset{(k_1)}{(1,1)},...,(1,1),\underset{(k_2)}{(0,1)},...,(0,1),
\underset{(k_3)}{(0,0)},...\big).
\end{align*}

Alors
\begin{align*}
h(\sigma_{\partial\delta})=\sigma_{\partial h(
\sigma_\delta)}+(-1)^{\ell_1+\dots+\ell_{k_2-1}}\left( \sum_{1\leq
i<k_0}\ell_i+\sum_{k_0\leq i<k_1}(\ell_i-f_i)\right)
\sigma_\delta.
\end{align*}
\end{prop}

\vskip 0.3cm
\noindent
{\bf Preuve}\\

On reprend les notations de la proposition pr\'ec\'edente. On a
donc:

\begin{align*}
\sigma_{\partial\delta}&=-\sum_{\Delta\in\sigma_\delta}a_{\Delta}
\Big\{\sum_{\begin{smallmatrix}i<j\cr0\leq i<k_0-1\cr k_1-1\leq
j\leq k_2-1\end{smallmatrix}}\sum_{\Delta_{i,j}^{I,\ell_{i+1}-1,s}
\begin{smallmatrix}prop\cr\rightarrow\cr i\end{smallmatrix}\Delta}
\varepsilon_{i,j}^{\ell_{i+1}-1,s}\Delta_{i,j}^{I,\ell_{i+1}-1,s}\cr
&+\sum_{\begin{smallmatrix}i<j\cr k_0-1\leq i<k_1-1\cr k_1-1\leq
j\leq k_2-1\end{smallmatrix}}\Big[\sum_{\Delta_{j,i}^{I,0,1}
\begin{smallmatrix}prop\cr\rightarrow\cr
i\end{smallmatrix}\Delta}\varepsilon_{j,i}^{0,1}\Delta_{j,i}^{I,0,1}+
\sum_{\Delta_{i,j}^{I,\ell_{i+1}-1,s}\begin{smallmatrix}prop\cr\rightarrow\cr
i
\end{smallmatrix}\Delta}\varepsilon_{i,j}^{\ell_{i+1}-1,s}
\Delta_{i,j}^{I,\ell_{i+1}-1,s}\Big]\cr
&+\sum_{\begin{smallmatrix}i<j\cr k_1-1\leq i<k_2-1\cr k_1-1\leq
j\leq k_2-1
\end{smallmatrix}}\Big[\sum_{\Delta_{j,i}^{I,0,1}\begin{smallmatrix}prop\cr
\rightarrow\cr
i\end{smallmatrix}\Delta}\varepsilon_{j,i}^{0,1}\Delta_{j,i}^{I,0,1}
+\sum_{\Delta_{i,j}^{I,0,1}\begin{smallmatrix}prop\cr\rightarrow\cr
i\end{smallmatrix}\Delta}\varepsilon_{i,j}^{0,1}\Delta_{i,j}^{I,0,1}\Big]\Big\}.
\end{align*}

Le dernier sommet d'ordre $(1,1)$ dans $\partial\delta$ est
$k_2-1$. Dans l'expression pr\'ec\'edente, on regarde les termes
pour $j=k_2-1$:

\begin{align*}
&=-\sum_{\Delta\in\sigma_\Delta}
a_\Delta\Big\{\sum_{\begin{smallmatrix}i<j\cr0\leq i<k_0-1\cr
j=k_2-1
\end{smallmatrix}}\sum_{\Delta_{i,k_2-1}^{I,\ell_{i+1}-1,s}\begin{smallmatrix}prop
\cr\rightarrow\cr
i\end{smallmatrix}\Delta}\varepsilon_{i,k_2-1}^{\ell_{i+1}-1,s}
\Delta_{i,k_2-1}^{I,\ell_{i+1}-1,s}\cr
&+\sum_{\begin{smallmatrix}i<j\cr k_0-1\leq i<k_1-1\cr
j=k_2-1\end{smallmatrix}}
\Big[\sum_{\Delta_{k_2-1,i}^{I,0,1}\begin{smallmatrix}prop\cr\rightarrow\cr
i
\end{smallmatrix}\Delta}\varepsilon_{k_2-1,i}^{0,1}\Delta_{k_2-1,i}^{I,0,1}+
\sum_{\Delta_{i,k_2-1}^{I,\ell_{i+1}-1,s}\begin{smallmatrix}prop\cr\rightarrow\cr
i\end{smallmatrix}\Delta}\varepsilon_{i,k_2-1}^{\ell_{i+1}-1,s}
\Delta_{i,k_2-1}^{I,\ell_{i+1}-1,s}\Big]\cr
&+\sum_{\begin{smallmatrix}i<j\cr k_1-1\leq i<k_2-1\cr j= k_2-1
\end{smallmatrix}}\Big[\sum_{\Delta_{k_2-1,i}^{I,0,1}\begin{smallmatrix}prop\cr
\rightarrow\cr i\end{smallmatrix}\Delta}
\varepsilon_{k_2-1,i}^{0,1}\Delta_{k_2-1,i}^{I,0,1}+\sum_{
\Delta_{i,k_2-1}^{I,0,1}\begin{smallmatrix}prop\cr\rightarrow\cr
i\end{smallmatrix}\Delta}\varepsilon_{i,k_2-1}^{0,1}\Delta_{i,k_2-1}^{I,0,1}\Big]
\Big\}.
\end{align*}

{\bf *} Pour $0\leq i<k_0-1$,

\begin{align*}
\Delta_{i,k_2-1}^{I,\ell_{i+1}-1,s}=&\Big(\Big(
\overrightarrow{0c_1^0}\dots\overrightarrow{0c^0_{\ell_1}}\Big),\dots,
\Big(\overrightarrow{ic_1^i}\dots\underbrace{\overrightarrow{i(k_2-1)}}_{(s)}
\dots\overrightarrow{ic^i_{\ell_{i+1}-1}}\Big),
\Big(\overrightarrow{(k_2-1)c^i_{\ell_{i+1}}}\Big),\dots\cr
&\dots\widehat{k_2-1}\dots,
\Big(\overrightarrow{nc_1^n}\dots\overrightarrow{nc^n_{\ell_n}}\Big)\Big).
\end{align*}

Appliquons l'homotopie $h$, on obtient

\begin{align*}
&h(\Delta_{i,k_2-1}^{I,\ell_{i+1}-1,s})=\cr
&\hskip-0.2cm\left(\left(\overrightarrow{0c_1^0}\dots
\overrightarrow{0c^0_{\ell_1}}
\right),\dots,\left(\overrightarrow{ic_1^i}\dots
\underbrace{\overrightarrow{i(c^i_{\ell_{i+1}})}}_{(s)}\dots
\overrightarrow{ic^i_{\ell_{i+1}-1}}
\right),\dots\widehat{k_2-1}\dots,\left(\overrightarrow{nc_1^n}\dots
\overrightarrow{nc^n_{\ell_n}} \right)\right)
\end{align*}
et le signe sera

$$
\varepsilon_{i,k_2-1}^{\ell_{i+1}-1,s}=(-1)^{s-1}(-1)^{\ell_{i+1}}(-1)^{\ell_1+
\dots+\ell_{k_2-1}}.
$$
D'o\`u, pour chaque $1\leq s\leq \ell_{i+1}$,

$$
\varepsilon_{i,k_2-1}^{\ell_{i+1}-1,s}h(\Delta_{i,k_2-1}^{I,\ell_{i+1}-1,s})=
-(-1)^{\ell_1+\dots+\ell_{k_2-1}}\Delta(0,\dots\widehat{k_2-1}\dots,n).
$$\\

{\bf *} Pour $k_0-1\leq i<k_1-1$.

- On a d'abord $f_{i+1}$ graphes de type $\Delta_{k_2-1,i}$

\begin{align*}
&\Delta_{k_2-1,i}^{I,0,1}=\cr
&\left(\left(\overrightarrow{0c_1^0}\dots
\overrightarrow{0c^0_{\ell_1}}
\right),\dots,\left(\overrightarrow{(k_2-1)i}\right),\left(
\overrightarrow{ic_{1}^i}\dots\overrightarrow{ic^i_{\ell_{i+1}}}
\right),\dots\widehat{k_2-1}\dots,\left(\overrightarrow{nc_1^n}\dots
\overrightarrow{nc^n_{\ell_n}} \right)\right)\end{align*} alors
\begin{align*}
&h(\Delta_{k_2-1,i}^{I,0,1})=\cr&\left(\left(\overrightarrow{0c_1^0}\dots
\overrightarrow{0c^0_{\ell_1}}\right),\dots,\left(\overrightarrow{ic_{1}^i}
\dots\overrightarrow{ic^i_{\ell_{i+1}}}\right),\dots\widehat{k_2-1}\dots,
\left(\overrightarrow{nc_1^n}\dots\overrightarrow{nc^n_{\ell_n}}\right)\right).
\end{align*}
D'o\`u,

$$
\varepsilon_{k_2-1,i}^{0,1}h(\Delta_{k_2-1,i}^{I,0,1})=(-1)^{\ell_1+\dots+
\ell_{k_2-1}}\Delta(0,\dots\widehat{k_2-1}\dots,n).
$$\\

- On a aussi un graphe de type $\Delta_{i,k_2-1}$
\begin{align*}
\Delta_{i,k_2-1}^{I,\ell_{i+1}-1,s}=&\Big(\Big(
\overrightarrow{0c_1^0}\dots\overrightarrow{0c^0_{\ell_1}}\Big),\dots,
\Big(\overrightarrow{ic_1^i}\dots\underbrace{\overrightarrow{i(k_2-1)}}_{(s)}
\dots\overrightarrow{ic^i_{\ell_{i+1}-1}}\Big),
\Big(\overrightarrow{(k_2-1)c^i_{\ell_{i+1}}}\Big),\dots\cr
&\dots\widehat{k_2-1}\dots,
\Big(\overrightarrow{nc_1^n}\dots\overrightarrow{nc^n_{\ell_n}}\Big)\Big).
\end{align*}
D'o\`u, pour chaque $1\leq s\leq \ell_{i+1}$,

$$
\varepsilon_{i,k_2-1}^{\ell_{i+1}-1,s}h(\Delta_{i,k_2-1}^{I,\ell_{i+1}-1,s})=
-(-1)^{\ell_1+\dots+\ell_{k_2-1}}\Delta(0,\dots\widehat{k_2-1}\dots,n).
$$\\

{\bf *} Pour $k_1-1\leq i<k_2-1$, on a aussi des graphes des deux types :

\begin{align*}
\Delta_{k_2-1,i}^{I,0,1}=\left(\left(\overrightarrow{0c_1^0}\dots
\overrightarrow{0c^0_{\ell_1}}\right),\dots,\left(\overrightarrow{(k_2-1)i}\right),
\left(\overrightarrow{ic_{1}^i}\right),\dots\widehat{k_2-1}\dots,\left(
\overrightarrow{nc_1^n}\dots\overrightarrow{nc^n_{\ell_n}}\right)\right).
\end{align*}
Donc

$$
\varepsilon_{k_2-1,i}^{0,1}h(\Delta_{k_2-1,i}^{I,0,1})=(-1)^{\ell_1+\dots+
\ell_{k_2-1}}\Delta(0,\dots\widehat{k_2-1}\dots,n).
$$\\

De m\^eme, on a

\begin{align*}
\Delta_{i,k_2-1}^{I,0,1}=\left(\left(\overrightarrow{0c_1^0}\dots
\overrightarrow{0c^0_{\ell_1}}\right),\dots,\left(\overrightarrow{i(k_2-1)}
\right),\left(\overrightarrow{(k_2-1)c_{1}^i}\right),\dots\widehat{k_2-1}\dots,
\left(\overrightarrow{nc_1^n}\dots\overrightarrow{nc^n_{\ell_n}}\right)\right).
\end{align*}
Donc

$$
\varepsilon_{i,k_2-1}^{0,1}h(\Delta_{i,k_2-1}^{I,0,1})=-(-1)^{\ell_1+\dots+
\ell_{k_2-1}}\Delta(0,\dots\widehat{k_2-1}\dots,n).
$$\\

Enfin, dans l'expression de $\sigma_{\partial\Delta}$, on regarde
les termes pour $j<k_2-1$. D'abord, par d\'efinition de
l'homotopie $h$ si $j<k_2-1$, on a:
$$
\Delta_{i,j}^{I,r,s}\begin{matrix}prop\cr\rightarrow\cr
i\end{matrix}\Delta~~
\Longleftrightarrow~~h(\Delta_{i,j}^{I,r,s})\begin{matrix}prop\cr\rightarrow\cr
i
\end{matrix}h(\Delta)
$$
et
$$
\Delta_{j,i}^{I,r,t}\begin{matrix}prop\cr\rightarrow\cr
i\end{matrix}\Delta~~
\Longleftrightarrow~~h(\Delta_{j,i}^{I,r,t})\begin{matrix}prop\cr\rightarrow\cr
i
\end{matrix}h(\Delta).
$$
De plus les signes $\varepsilon^{r,s}_{i,j}$ et
$\varepsilon^{r,t}_{j,i}$ ne d\'ependent pas de $k_2-1$, on en
d\'eduit que la somme de tous les termes pour $j<k_2-1$ est:

\begin{align*}
\sigma_{\partial\delta}&=-\sum_{\Delta\in\sigma_\delta}a_{\Delta}
\Big\{\sum_{\begin{smallmatrix}i<j\cr0\leq i<k_0-1\cr k_1-1\leq
j<k_2-1\end{smallmatrix}}\sum_{\Delta_{i,j}^{I,\ell_{i+1}-1,s}
\begin{smallmatrix}prop\cr\rightarrow\cr i\end{smallmatrix}\Delta}
\varepsilon_{i,j}^{\ell_{i+1}-1,s}\Delta_{i,j}^{I,\ell_{i+1}-1,s}\cr
&+\sum_{\begin{smallmatrix}i<j\cr k_0-1\leq i<k_1-1\cr k_1-1\leq
j< k_2-1\end{smallmatrix}}\Big[\sum_{\Delta_{j,i}^{I,0,1}
\begin{smallmatrix}prop\cr\rightarrow\cr
i\end{smallmatrix}\Delta}\varepsilon_{j,i}^{0,1}\Delta_{j,i}^{I,0,1}+
\sum_{\Delta_{i,j}^{I,\ell_{i+1}-1,s}\begin{smallmatrix}prop\cr\rightarrow\cr
i
\end{smallmatrix}\Delta}\varepsilon_{i,j}^{\ell_{i+1}-1,s}
\Delta_{i,j}^{I,\ell_{i+1}-1,s}\Big]\cr
&+\sum_{\begin{smallmatrix}i<j\cr k_1-1\leq i<k_2-1\cr k_1-1\leq
j< k_2-1
\end{smallmatrix}}\Big[\sum_{\Delta_{j,i}^{I,0,1}\begin{smallmatrix}prop\cr
\rightarrow\cr
i\end{smallmatrix}\Delta}\varepsilon_{j,i}^{0,1}\Delta_{j,i}^{I,0,1}
+\sum_{\Delta_{i,j}^{I,0,1}\begin{smallmatrix}prop\cr\rightarrow\cr
i\end{smallmatrix}\Delta}\varepsilon_{i,j}^{0,1}\Delta_{i,j}^{I,0,1}\Big]\Big\}\cr
&=\sigma_{\partial
h(\sigma_\Delta)}(0,\dots\widehat{k_2-1}\dots,n).
\end{align*}

Ceci d\'emontre notre proposition.


\

\section{Cohomologie des graphes ascendants}\label{sec6}

\

Disons qu'un sommet $i$ d'un graphe de type $(\pt\leq\pt)$
$\Delta$ est \underbar{simple} si $\#Deb(i)\leq1$ et
$\#Fin(i)\leq1$.

\vskip 0.3cm
\begin{prop}\label{simple} {\rm (Modulo un cobord, le symbole d'un cocycle ne
contient que des sommets simples)}

\

Soit $\delta$ une combinaison lin\'eaire sym\'etrique de graphes
de type $(\pt\leq\pt)$. On suppose que $\delta$ est un cocycle
($\partial\delta=0$). Alors il existe un cobord $\partial\beta$
tel que le symbole $\sigma_{\delta-\partial\beta}$ de
$\delta-\partial\beta$ ne contient que des graphes $\Delta$ dont
tous les sommets sont simples.
\end{prop}

\vskip 0.3cm
\noindent
{\bf Preuve}

Puisque $\delta$ est un cocycle, $h(\sigma_{\partial\delta})=0$.
On a alors:
$$
0=(-1)^{\ell_1+\dots+\ell_{k_2-1}}\left(\sum_{1\leq i<k_0}\ell_i+
\sum_{k_0\leq
i<k_1}(\ell_i-f_i)\right)\sigma_\delta+\sigma_{\partial
h(\sigma_\delta)}=a\sigma_\delta+\sigma_{\partial
h(\sigma_\delta)}.
$$
Dire que $\sigma_\delta$ contient au moins un graphe poss\'edant
un sommet non simple, c'est dire que dans ${\mathcal O}(\delta)$,
$k_1>0$. Dans l'expression pr\'ec\'edente, le coefficient $a$ de
$\sigma_\delta$ est non nul. Puisque $\sigma_\delta$ n'est pas
nul, on en d\'eduit que $h(\sigma_\delta)$ n'est pas nul, donc
$k_2>k_1$, il existe des sommets d'ordre $(1,1)$ dans les graphes
$\Delta$ de $\sigma_\delta$.

Posons $\beta_1=-\displaystyle\frac{1}{a}S(h(\sigma_\delta))$, où
$S(h(\sigma_\delta))$ est la symétrisation du graphe
$h(\sigma_\delta)$. D'abord par construction, les graphes
apparaissant dans $S(h(\sigma_\delta))$ mais pas dans
$h(\sigma_\delta)$ sont d'ordre strictement plus petit que les
graphes apparaissant dans $h(\sigma_\delta)$, ou:
$$
\sigma_{h(\sigma_\delta)}=\sigma_{S(h(\sigma_\delta))}.
$$
Ensuite, on a vu que pour calculer le symbole de $\partial\delta$,
on ne consid\'erait que les graphes du symbole de $\delta$ donc:
$$
\sigma_{\partial h(\sigma_\delta)}=\sigma_{\partial S(h(\sigma_\delta))}=-
a\sigma_{\partial\beta_1}
$$

Alors,
$$
\sigma_{\partial\beta_1}=\sigma_\delta.
$$
Autrement dit $\delta_1=\delta-\partial\beta_1$ a un ordre strictement plus
petit que ${\mathcal O}(\delta)$. Si le symbole de $\delta_1$ a des graphes avec
des sommets non simples, on peut recommencer cette op\'eration. Au bout d'un
nombre fini d'\'etapes, on arrive sur une combinaison de graphes $\delta-
\partial\beta$ dont le symbole ne contient que des sommets simples:
$$
{\mathcal
O}(\delta-\partial\beta)=\Big((1,1),\dots,(1,1),(0,1),\dots,(0,1),(0,0),\dots,(0,0)\Big).
$$
Donc tous les graphes de $\delta-\partial\beta$ n'ont eux aussi
que des sommets simples. Ces graphes sont des graphes vectoriels
ou linéaires \`a sommets simples \'etudi\'es dans \cite{[AAC1]} et
\cite{[AAC2]} .

\vskip 0.3cm
\begin{defn} {\rm (Les roues)}

\

Une roue sym\'etrique $R_k$ de longueur $k$ est le sym\'etris\'e
de la roue simple qui est le graphe $\Delta_k$ ayant $k$ sommets
$(1,\dots,k)$ et les $k$ fl\`eches
$\Big((\overrightarrow{12}),(\overrightarrow{23}),\dots,
(\overrightarrow{(k-1)k}),(\overrightarrow{k1})\Big)$.
\end{defn}

\vskip  0.3cm
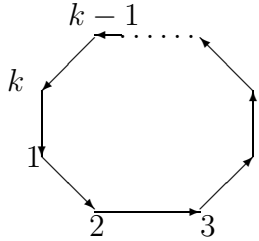
\begin{center}
\begin{picture}(200,100)(100,50)
\put(200,146){\vector(-1,0){10}}
\put(190,145){\vector(-1,-1){20}} \put(170,125){\vector(0,-1){25}}
\put(190,79){\vector(1,0){40}} \put(170,100){\vector(1,-1){20}}
\put(250,100){\vector(0,1){25}}
\put(250,125){\vector(-1,1){20}}\put(199,145){\dots\dots}
\put(230,80){\vector(1,1){20}}\put(164,96){$1$}
\put(188,70){$2$}\put(157,125){$k$}\put(180,150){$k-1$}\put(230,70){$3$}
\end{picture}
\end{center}

\

Le calcul de la cohomologie des graphes vectoriels et linéaires a
\'et\'e r\'ealis\'e dans \cite{[AAC1]} et \cite{[AAC2]}. Comme
nous nous sommes ramen\'es \`a cette situation, on obtient, de
m\^eme, ici:

\vskip 0.3cm
\begin{thm}{\rm (La cohomologie de Chevalley des graphes de type $(\pt\leq\pt)$)}

\

La cohomologie de Chevalley des graphes de type $(\pt\leq\pt)$ est
donn\'ee par les roues de longueur impaire. Plus pr\'ecis\'ement,
pour tout $n$, une base de $H^n$ est donn\'ee par
$$
\left\{R_{2k_1+1}\wedge R_{2k_2+1}\wedge\dots\wedge R_{2k_p+1},\quad\hbox{avec}
~~k_1<k_2<\dots<k_p,~~\sum_{i=1}^p(2k_i+1)=n\right\}.
$$
\end{thm}

\vskip 0.3cm

\

\begin{thm}{\rm (La cohomologie de Chevalley de l'action adjointe des tenseurs
ascendants)}

\

La cohomologie de Chevalley de l'action adjointe des tenseurs
ascendants dans $T_{poly}(\mathbb{R}^d)$ est donn\'ee par des
cocha\^ines $C_\delta$ o\`u $\delta$ est une combinaison de roues
impaires de longueur strictement plus petite que $2d$.

\end{thm}

\vskip 0.3cm

\noindent {\bf Preuve}

Soit $C_\delta$ un cocycle d\'efini avec des graphes de type
$(\pt\leq\pt)$, alors, $\partial C_\delta=C_{\partial\delta}=0$.
D'apr\`es le th\'eor\`eme pr\'ec\'edent,
$$
\delta=\partial\beta+\sum_{r_1<\dots<r_q\atop r_j impair}a_{r_1\dots
r_q}R_{r_1}\wedge\dots \wedge R_{r_q}.
$$

Montrons que si $2p+1>2d$, alors, $C_{R_{2p+1}}=0$.

Prenons $\Delta$ la roue de longueur $2p+1$:
$$
\Delta=\left(\left(\overrightarrow{12}\right),\left(\overrightarrow{23}\right)
,\dots\left(\overrightarrow{(2p+1)1}\right)\right).
$$
Alors si $\alpha_1,\dots,\alpha_{2p+1}$ sont des tenseurs d'ordre
$k_1,\dots,k_{2p+1}$. Posons

\noindent
$I_1=\{1,\dots,k_1\},I_2=\{k_1+1,\dots,k_1+k_2\},\dots,I_{2p+1}=\{\displaystyle
\sum_{j<2p+1}k_j+1, \dots,\sum_{j\leq2p+1}k_j\}$. Soient
$(j_1,\dots,j_{2p+1})\in I_1\times\dots\times I_{2p+1}$, posons
$\widehat{I_1}=I_1\setminus
\{j_1\},\dots,\widehat{I_{2p+1}}=I_{2p+1}\setminus\{j_{2p+1}\}.$
Alors,
\begin{align*}
&C_\Delta(\alpha_1,\dots,\alpha_{2p+1})=\cr&\sum_{1\leq t_1,\dots,
t_{k_1+\dots+k_{2p+1}}\leq d }~~\sum_{(j_1,\dots,j_{2p+1})\in
I_1\times\dots\times
I_{2p+1}}\varepsilon(j_1,\dots,j_{2p+1},\widehat{I_1},\dots,\widehat{I_{2p+1}})
\cr&\partial_{t_{j_{2p+1}}}
\alpha_1^{t_{I_1}}\partial_{t_{j_1}}\alpha_2^{t_{I_2}}\dots
\partial_{t_{j_{2p}}}\alpha_{2p+1}^{t_{I_{2p+1}}}\partial_{t_{\widehat{I_1}}}
\wedge\dots\wedge\partial_{t_{\widehat{I_{2p+1}}}}.
\end{align*}
On d\'efinit $2p+1$ champs de vecteurs d\'ependant de
$t_{\widehat{I_1}},\dots,t_{\widehat{I_{2p+1}}}$ en posant
$$
\xi_1=\sum_{t_{j_1}=1}^d\alpha_1^{t_{I_1}}\partial_{t_{j_1}}~~,\dots,~~\xi_{2p+1}
=\sum_{t_{j_{2p+1}}=1}^d\alpha_{2p+1}^{t_{I_{2p+1}}}\partial_{t_{j_{2p+1}}}
$$
et
$$
Jac_{j_1,\dots,j_{2p+1}}(\xi_k)=\begin{pmatrix}_{\partial_{t_{j_1}}\xi_k^{t_{j_1}}}
&_{.}&_{.}&_{.}&_{.}&_{\partial_{t_{j_{2p+1}}}\xi_k^{t_{j_1}}}\\_{\partial_{t_{j_1}}
\xi_k^{t_{j_2}}}& _{.} & _{.} &_{.}&_{.}&_{\partial_{t_{j_{2p+1}}}
\xi_k^{t_{j_2}}} \\_{.} & _{} & _{} & _{} & _{} & _{.} \\
_{.}&_{}&_{}& _{} & _{}& _{.} \\ _{.} & _{} & _{} & _{} & _{} &
_{.} \\ _{\partial_{t_{j_1}}\xi_k^{t_{j_{2p+1}}}} & _{.} & _{.} &
_{.}&_{.}&_{\partial_{t_{j_{2p+1}}}\xi_k^{t_{j_{2p+1}}}}\end{pmatrix}.
$$

En sym\'etrisant $C_\Delta$, on obtient

\begin{align*}
&C_{R_{2p+1}}(\alpha_1,\dots,\alpha_{2p+1})=\sum_{\sigma\in
S_{2p+1}}\varepsilon_{|\alpha|}(\sigma)C_\Delta(\alpha_{\sigma(1)},
\dots,\alpha_{\sigma(2p+1)})\cr&=\sum_{1\leq t_1,\dots,
t_{k_1+\dots+k_{2p+1}}\leq d }~~\sum_{(j_1,\dots,j_{2p+1})\in
I_1\times\dots\times
I_{2p+1}}\varepsilon(j_1,\dots,j_{2p+1},\widehat{I_1},\dots,\widehat{I_{2p+1}})
\cr&\hskip 1.5cm{\mathfrak
a}\Big(tr\big(Jac_{j_1,\dots,j_{2p+1}}(\xi_1)\dots
Jac_{j_1,\dots,j_{2p+1}}(\xi_{2p+1})\big)\Big)\partial_{t_{\widehat{I_1}}}
\wedge\dots\wedge\partial_{t_{\widehat{I_{2p+1}}}}
\end{align*}

o\`u on a not\'e par $\mathfrak{a}$ l'op\'erateur
d'antisym\'etrisation.

\vskip 0.3cm

Dans \cite{[W]} et \cite{[AL]}, il est montr\'e que
$\mathfrak{a}(tr(A_1\dots A_{2p+1}))=0$ d\`es que $2p+1>2d$.
Alors, ceci prouve que $C_{R_{2p+1}}=0$ si $2p+1>2d.$

\

 La cohomologie est donc engendr\'ee par les cocha\^\i nes $C_\delta$
 o\`u $\delta$ est une combinaison de roues
impaires de longueur strictement plus petite que $2d$.

\vskip 0.3cm

\begin{rema}
\

On peut traiter de m\^eme les graphes de Kontsevich aériens de
type $(\pt\geq\pt)$ formés par des sommets de type
$(f_i\geq\ell_i)$, c'est \`a dire, exactement $\ell_i$ fl\`eches
sortent de $i$ et $f_i$ fl\`eches aboutissent sur $i$ avec
$\begin{cases}
  f_i>\ell_i & \text{si $f_i>1$}, \\
    f_i\geq \ell_i & \text{si $f_i\leq1$}.
  \end{cases}$.
\vskip0.2cm

On peut regarder la cohomologie de Chevalley des graphes de type
$(\pt\geq\pt)$ en d\'efinissant un autre ordre:\vskip0.2cm

Classe 1: les sommets $i$ de type $(f_i\geq0)$ avec $f_i>1$. On
posera $\mathcal{O}(i)=(f_i,0)$.\\

Classe 2: les sommets $i$ de type $(f_i\geq\ell_i)$ avec
$f_i>\ell_i\geq1$. On
posera $\mathcal{O}(i)=(f_i,\ell_i)$.\\

Classe 3: les sommets $i$ de type $(1\geq1)$. On
posera $\mathcal{O}(i)=(1,1)$.\\

Classe 4: les sommets $i$ de type $(1\geq0)$. On
posera $\mathcal{O}(i)=(1,0)$.\\

Classe 5: les sommets $i$ de type $(0\geq0)$. On
posera $\mathcal{O}(i)=(0,0)$.\\

On ordonne les ordres des sommets suivant $(*)$ en posant:
\begin{align*}
&(f_i,0)>(f_j,\ell_j)>(1,1)>(1,0)>(0,0), \\& ~~(f_i,0)\geq
(f_{i'},0)\Leftrightarrow f_i\geq
f_{i'}~~\hbox{et}~~(f_i,\ell_i)\geq
(f_{i'},\ell_{i'})\Leftrightarrow
  \begin{cases}
    f_i\geq f_{i'} , \\
    f_i= f_{i'}~~ \hbox{et}~~ \ell_i\geq \ell_{i'}.
  \end{cases}
\end{align*}

L'ordre ${\mathcal O}(\Delta)$ d'un graphe  $\Delta$ de type
$(\pt\geq\pt)$ est le mot form\'e par les ordres de ses sommets:
$$
{\mathcal O}(\Delta)=\left({\mathcal O}(1),\dots,{\mathcal
O}(n)\right).
$$
On ordonne les ordres des graphes par l'ordre lexicographique, en
respectant $(*)$.

En appliquant le m\^eme technique, on retrouve que la cohomologie
des graphes de type $(\pt\geq\pt)$ est aussi engendr\'ee par une
combinaison de roues impaires.

\

On remarque, finalement, qu'en consid\'erant les graphes de type
$(\pt\leq\pt)$ et de type $(\pt\geq\pt)$, on passe d'une classe de
graphes \`a l'autre en changeant le sens de toutes les fl\`eches,
on dira qu'on a transpos\'e ce graphe. Les roues sont des graphes
invariants sous cette transformation, elles donnent, aussi, la
base de la cohomologie transpos\'ee d\'efinie par:
$$
(~~^t\partial)(\Delta)=~~^t(\partial(~~^t\Delta)).
$$
Cependant, l'op\'erateur $~~^t\partial$ est un op\'erateur de
cobord diff\'erent de celui de Chevalley sur $T_{poly}({\mathbb
R}^d)$. (Voir l'exemple donné dans \cite{[AAC2]})

\end{rema}


\

\end{document}